\documentclass[11pt]{amsart}

\usepackage[active]{srcltx} 
\usepackage{amsmath}
\usepackage{amsfonts}
\usepackage{amssymb}
\usepackage{latexsym}
\usepackage{graphicx}
\usepackage[cp866]{inputenc}

 \advance\hoffset by -0.0 truecm

\newtheorem{theorem}{Theorem}[section]

\newtheorem{lemma}[subsection]{Lemma}

\theoremstyle{definition}
\newtheorem{definition}[subsection]{Definition}
\theoremstyle{remark}
\newtheorem{remark}[theorem]{Remark}
\numberwithin{equation}{section}

\newcommand{\R}{\mathbb{R}}
\newcommand{\C}{\mathbb{C}}
\newcommand{\ii}{\mathbf{i}}

\newcommand{\ub}{\mathbf u}
\newcommand{\vb}{\mathbf v}

\DeclareMathOperator{\re}{\hbox{Re}}
\DeclareMathOperator{\im}{\hbox{Im}}

\begin{document}
\title[Quaternions operator $\Delta_{\lambda}$]
{Quaternion $H$-type group and differential operator
$\Delta_{\lambda}$}

\author{Der-Chen Chang,\ Irina Markina}

\address{Department of Mathematics, Georgetown University, Washington
D.C. 20057, USA}

\email{chang@georgetown.edu}

\address{Department of Mathematics,
University of Bergen, Johannes Brunsgate 12, Bergen 5008, Norway}

\email{irina.markina@uib.no}

\thanks{The first author is partially supported by
a Competitive Research Grant at Georgetown University. The second
author is partially supported by a Research Grant at University of
Bergen.}




\begin{abstract}

We study the relations between the quaternion $H$-type group and the
boundary of the unit ball on two dimensional quaternionic space. The
orthogonal projection of the space of square integrable functions
defined on quaternion $H$-type group into its subspace of boundary
values of $q$-holomorphic functions is consider. The precise form of
Cauchy-Szego kernel and the orthogonal projection operator is
obtained. The fundamental solution for the operator
$\Delta_{\lambda}$ is found.

\end{abstract}

\maketitle

\centerline {\small \small \it Dedicated to Professor Qikeng Lu on
his 80th birthday}
\section{Introduction}

The real number system is extended to the complex then to the
quaternion systems of numbers and further finds its the most
exciting generalization  incorporated the geometric concept of the
direction, the so-called Clifford algebras, which Clifford himself
called ``geometric algebras",~\cite{Loun,Por}, see also~\cite{GS}
for Clifford analysis and numerous  references. The division
algebras, which are only real, complex, quaternionic and
octonionic numbers give the origin to homogeneous groups
satisfying $J^2$ condition~\cite{CDKR}.  The simplest
non-commutative example of them is the Heisenberg group closely
related to the complex number system and that found its numerous
applications in physics, quantum mechanics, differential geometry
(see, for instance~\cite{CChGr3,Kap,Kor}). The following more
complicate example is an quaternion analogue of the Heisenberg
group that has at least four dimensional horizontal distribution
and three dimensional center~\cite{ChM}. We call this group
quaternion $H$-type group due to~\cite{Kap1}.

The quaternion $H$-type group $\mathcal Q$ can be realized as a
boundary of the unit ball on the two dimensional quaternionic
space $\mathbb H^2$. The Siegel upper half space $\mathcal U_1\in\mathbb H^2$ is
$q$-holomorphically  equivalent to the unit ball in $\mathbb H^2$.
The group $\mathcal Q$ arise as the group of translations of
$\mathcal U_1$. This leads to its identification with the boundary
$\partial \mathcal U_1$. We give the precise formulas of the
action. Because of this identification and by the use of various
symmetries of $\mathcal U_1$ the Cauchy-Szeg\"{o} projector is
realized as a convolution operator on the group $\mathcal Q$ with
explicitly given singular kernel. The analogue for the group
$\mathcal Q$ of the Laplace operator is the so called
$\Delta_{\lambda}$ operator that is expressed as a sum of the
square of vector fields forming the frame of the horizontal
distribution plus the $\lambda$ times their commutators. We find
the kernel of the operator $\Delta_{\lambda}$.

Part of this article is based on a lecture presented by the first
author during the International Conference on Several Complex
Variables which was held on June 5-9, 2006 at the Chinese Academy of
Sciences, Beijing, China. The first author thanks the organizing
committee, especially Professor Xiangyu Zhou for his invitation. He
would also like to thank all the colleagues at the Institute of
Mathematics, AMSS, Chinese Academy of Sciences for the warm
hospitality during his visit to China. We also would like to thank
Professor Jingzhi Tie for many inspired conversations on this
project.

\section{Quaternion $H$-type group and the Siegel upper half space.}

We remember shortly the definitions of quaternion numbers. Let
$i_1,i_2,i_3$ be three imaginary units such that
$$i^2_1=i_2^2=i_3^2=i_1i_2i_3=-1.$$ The multiplication between the imaginary units is given in Table~\ref{t1}.
\begin{table}[ht]
\caption{Multiplication of imaginary units} \label{t1}
\begin{center}
\begin{tabular}{|| c | r | r | r ||}
\hline \     &  $i_1$ & $i_2$ & $i_3$
\\  \hline
$i_1$ &  $-1$ & $i_3$ & $-i_2$
\\ \hline
$i_2$ & $-i_3$ & $-1$ & $i_1$
\\ \hline
$i_3$ & $i_2$ & $-i_1$ & $-1$
\\ \hline
\end{tabular}
\end{center}
\end{table}
Any quaternion $q$ can be written in the algebraic form as
$q=t+ai_1+bi_2+ci_3$, where $t,a,b,c$ are real numbers.  The number
$t$ is called the {\it real part} and denoted by $t=\re q$. The
vector $\mathbf u=(a,b,c)$ is the {\it imaginary part} of $q$. We
use the notations
$$a=\im_1 q,\quad b=\im_2 q,\quad c=\im_3 q,\quad\text{and}\quad \im q=\mathbf u=(a,b,c).$$
Similarly to complex numbers, vectors, and matrices, the addition of
two quaternions is equivalent to summing up the elements. Set
$q=t+\ub$, and $h=s+xi_1+yi_2+zi_3=s+\vb$. Then
$$q+h=(t+s)+(\ub+\vb)=(t+s)+(a+x)i_1+(b+y)i_2+(c+z)i_3.$$
Addition satisfies all the commutation and association rules of real
and complex numbers. The quaternion multiplication (the Grassmanian
product) is defined by
$$qh=(ts-\ub\cdot\vb)+(t\vb+s\ub+\ub\times\vb),$$
where $\ub\cdot\vb$ is the scalar product and $\ub\times\vb$ is the
vector product of $\ub$ and $\vb$. The multiplication is not
commutative because of the non-commutative vector product. The
conjugate $\bar q$ to $q$ is defined in a similar way as for the
complex numbers: $\bar q=t-ai_1-bi_2-ci_3$. Then the modulus $|q|$
of $q$ is given by $|q|^2=q\bar q=t^2+a^2+b^2+c^2$. The scalar
product between $q$ and $h$ is\begin{equation} \label{sp} \langle
q,h\rangle=\re(q\bar h)=ts+ax+by+cz.
\end{equation} We also have \begin{equation*}
\overline{qh}  =  \bar h\bar q,\quad
 |qh|  =  |q||h|,\quad
 q^{-1}  =  \frac{\bar q}{|q|^2}.
\end{equation*}

The imaginary units have the representation by $(4\times 4)$ real
matrices:
\begin{equation*}\ii_1=\left[\array{rrrr}0 & -1 & 0 & 0
\\ 1 & 0 & 0 & 0
\\
0 & 0 & 0 & -1
\\ 0 & 0 & 1 & 0\endarray\right],\qquad
\ii_2=\left[\array{rrrr} 0 & 0 & -1 & 0
\\
0 & 0 & 0 & 1
\\
1 & 0 & 0 & 0
\\
0 & -1 & 0 & 0\endarray\right],\qquad \ii_3=\left[\array{rrrr} 0 & 0
& 0 & -1
\\
0 & 0 & -1 & 0
\\
0 & 1 & 0 & 0
\\
1 & 0 & 0 & 0\endarray\right].\end{equation*} Then a quaternion can
be written in the matrix form
$$Q=\left[\array{rrrr}
t & -a & -b & -c
\\
a & t & -c & b
\\
b & c & t & -a
\\
c & -b & a & t\endarray\right]=t\mathbf U+a\ii_1+b\ii_2+c\ii_3,$$
where $\mathbf U$ is the unit $(4\times 4)$ matrix. Notice that
\begin{itemize}
\item[1.]{$\det Q=|q|^4$,}
\item[2.]{$Q^{T}=-Q$ represents the conjugate quaternion $\bar q$,}
\item[3.]{$Q^{-1}=-\frac{1}{\sqrt{\det Q}}Q$,}
\item[4.]{$Q^{-1}$ represents the inverse quaternion
$q^{-1}$.}\end{itemize} The nice description of algebraic and
geometric properties of quaternion a reader can find in the original
book of W.~Hamilton~\cite{Ham}.

Let $\mathbb H^2$ be a two dimensional vector space of pairs
$h=(h_1,h_2)$ over the field of real numbers with the norm
$\|h\|^2=h^2_1+h^2_2$. We describe the Siegel upper half space in
$\mathbb H^2$ carrying out the counterpart with the two dimensional
complex space $\mathbb C^2$.

Let $D_1$ denotes the unit ball in $\mathbb C^2$:
$$D_1=\{(w_1,w_2)\in\mathbb C^2:\ |w_1|^2+|w_2|^2<1\}.$$ The set
$$U_1=\{(z_1,z_2)\in\mathbb C^2:\ \re z_2>|z_1|^2\}$$ is the Siegel
half space. The Caley transformation
\begin{equation*}
 w_1=\frac{2z_1}{1+z_2},\qquad    w_2=\frac{1-z_2}{1+z_2},
\end{equation*} and its inverse \begin{equation*}
 z_1=\frac{w_1}{1+w_2},\qquad  z_2=\frac{1-w_2}{1+w_2},
\end{equation*} show that the unit ball $D_1$ and Siegel half space $U_1$ are
biholomorphically equivalent.

Now, take the unit ball $B_1$ in $\mathbb H^2$
$$B_1=\{(h_1,h_2)\in{\mathbb H}^2:\ |h_1|^2+|h_2|^2<1\}$$ and a
Siegel half space in $\mathbb H^2$
$$\mathcal U_1=\{(q_1,q_2)\in{\mathbb H}^2:\ \re q_2>|q_1|^2\}.$$
The Caley transformation, mapping the unit ball $B_1$ to the Siegel
half space $\mathcal U_1$ and vice versa has the form:
\begin{equation*}
 h_1=q_1\big(1+(1+q_2)^{-1}(1-q_2)\big)=\frac{2q_1(1+\bar q_2)}{|1+q_2|^2},\quad
 h_2=(1+q_2)^{-1}(1-q_2)=\frac{(1+\bar q_2)(1-q_2)}{|1+q_2|^2},
\end{equation*} and the inverse transformation \begin{equation*}
 q_1=h_1(1+h_2)^{-1}=\frac{h_1(1+\bar h_2)}{|1+h_2|^2},
\quad q_2=(1-h_2)(1+h_2)^{-1}=\frac{(1-h_2)(1+\bar h_2)}{|1+h_2|^2}.
\end{equation*} Since the multiplication of quaternion
is not commutative, it is possible to define another Caley's
transformation, but the geometry will be the same. The boundary of
$\mathcal U_1$ is $$\partial\mathcal U_1=\{(q_1,q_2)\in{\mathbb
H}^2:\ \re q_2=|q_1|^2\}.$$

We mention here three automorphisms of the domain $\mathcal U_1$,
dilation, rotation and translation. Let $q=(q_1,q_2)\in\mathcal
U_1$. For each positive number $\delta$ we define a {\it dilation}
$\delta\circ q$ by $$\delta\circ q=\delta\circ(q_1,q_2)=(\delta
q_1,\delta^2 q_2).$$ The non-isotropy of the dilation comes from the
definition of $\mathcal U_1$. For each unitary linear transformation
$\mathcal R$ on $\mathbb H$ we define the {\it rotation} $\mathcal
R(q)$ on $\mathcal U_1$ by $$\mathcal R(q)=\mathcal
R(q_1,q_2)=(\mathcal R(q_1),q_2).$$ Both, the dilation and rotation
give $q$-holomorphic (the definition of $q$-holomorphic mapping see
below) self mappings of $\mathcal U_1$ and extend to mappings on the
boundary $\partial \mathcal U_1$. Before we describe a translation
on $\mathcal U_1$, we introduce the quaternionic $H$-type group
denoted by $\mathcal Q$. This group consists of the set $$\mathbb
H\times\mathbb R^3=\{[w,t]:\ w\in\mathbb H,
t=(t_1,t_2,t_3)\in\mathbb R^3\}$$ with the multiplication law
\begin{equation}\label{lm}
[w,t_1,t_2,t_3]\cdot[\omega,s_1,s_2,s_3]=[w+\omega,t_1+s_1-2\im_1\bar
\omega w,t_2+s_2-2\im_2\bar \omega w,t_3+s_3-2\im_3\bar \omega w].
\end{equation} The law~\eqref{lm} makes $\mathbb H\times\mathbb R^3$ into Lie group with the
neutral element $[0,0]$ and the inverse element $[w,t]^{-1}$ given by $[w,t_1,t_2,t_3]^{-1}=[-w,-t_1,-t_2,-t_3]$.

To each element $[w,t]$ of $\mathcal Q$ we associate the following
$q$-holomorphic affine self mapping of~$\mathcal U_1$.
\begin{equation} \label{acgr} [w,t_1,t_2,t_3]:
(q_1,q_2)\mapsto(q_1+w,q_2+|w|^2+2\bar w
q_1+i_1t_1+i_2t_2+i_3t_3).\end{equation}  This mapping preserve the
following "height" function \begin{equation} \label{r} r(q)=\re
q_2-|q_1|^2.
\end{equation} In fact, since $|q_1+w|^2=|q_1|^2+|w|^2+2\re \bar w q_1$, we obtain
$$\re(q_2+|w|^2+2\bar w q_1)-|q_1+w|^2=\re q_2-|q_1|^2.$$ Hence, the transformation~\eqref{acgr}
maps $\mathcal U_1$ to itself and preserves the boundary $\partial \mathcal U_1$.

The reader can check that the mapping~\eqref{acgr} defines an action
of the group $\mathcal Q$ on the space $ \mathcal U_1$. If one
composes the mappings~\eqref{acgr}, corresponding to elements
$[w,t]$, $[\omega,s]\in\mathcal Q$, the resulting transformation
will correspond to the element $[w,t]\cdot[\omega,s]$.
Thus,~\eqref{acgr} gives us a realization of  $\mathcal Q$ as a
group of affine $q$-holomorphic bijections of $\mathcal U_1$. We can
identify the elements of $\mathcal U_1$ with the boundary via its
action on the origin $$h(0)=[w,t]: (0,0)\mapsto
(w,|w|^2+i_1t_1+i_2t_2+i_3t_3),$$ where $h=[w,t]$.  Thus $$\mathcal
Q \ni [w,t_1,t_2,t_3]\mapsto
(w,|w|^2+i_1t_1+i_2t_2+i_3t_3)\in\partial \mathcal U_1.$$

We may use the following coordinates $(q_1,t,r)=(q_1,t_1,t_2,t_3,r)$
on $\mathcal U_1$: $$\mathcal U_1\ni
(q_1,q_2)=(q_1,t_1,t_2,t_3,r),$$ where $$t_1=\im_1 q_2,\quad
t_2=\im_2 q_2,\quad t_3=\im_3 q_2,\quad r=r(q_1,q_2)=\re q_2-|q_1
|^2. $$ If $\re q_2=|q_1 |^2$ we get the  coordinates on the
boundary $\partial\mathcal U_1$ of the Siegel half space
$$\partial\mathcal U_1\ni (q_1,q_2)=(q_1,t_1,t_2,t_3),$$ where $t_m$
are as above and $r=r(q_1,q_2)=0$.

\section{Tangential Cauchy-Riemann-Fueter operators}

Before going further, we collect here the necessary definitions
concerning the quaternion calculus. $Q$-holomorphic functions on
$\mathbb H$ were studied by Fueter and his
collaborators~\cite{Fuet1,Fuet2}. The reader can find the account of
the theory of $q$-holomorphic functions in~\cite{Deav,GS,Sud}.

We continue exploit the analogy with the complex variables theory.
Let $M$ be a manifold of dimension $2n$ with the complex structure
$I$. Such kind of manifolds is called {\it complex manifold}. Let
$f: M\to\mathbb C$ be a differentiable function defined on the
complex manifold $M$. Write $f=f_0+if_1$, where $f_0,f_1:
M\to\mathbb R$. Then $f$ is called holomorphic, if
\begin{equation}\label{cr}df_0+I(df_1)=0\end{equation} on $M$, where~\eqref{cr} is called the Cauchy-Riemann equation.
Let, now, $M$ be a manifold of dimention $4n$. A hypercomplex
structure on $M$ is a triple $(I_1,I_2,I_3)$ on $M$, where $I_k$ is
a complex structure on $M$ and $I_1I_2=I_3$. If $M$ has hypercomplex
structure, it is called the hypercomplex manifold. Let $f: M\to
\mathbb H$ be a smooth function defined on a hypercomplex manifold
$M$. Then $f=f_0+i_1f_1+i_2f_2+i_3f_3$, where $f_0,\ldots, f_3$ are
smooth functions.
 We define a $q$-holomorphic function on $M$ to be a smooth function $f: M\to \mathbb H$ for which
\begin{equation}\label{crf}
df=df_0+I_1(df_1)+I_2(df_2)+I_3(df_3)=0.\end{equation}
Equation~\eqref{crf} is the natural quaternionic analogue of the
Cauchy-Riemann equation~\eqref{cr} and is called the
Cauchy-Riemann-Fueter equation.

The equation~\eqref{crf} can be written in terms of the partial
derivatives with respect to the quaternionic variable. We introduce
the differential operators
\begin{equation}\begin{split} \frac{\partial_l f}{\partial \bar
q}=\bar{\partial}_l f & =  \frac{1}{2} \Big(\frac{\partial
f}{\partial x_0}+\sum_{m=1}^{3}i_m\frac{\partial f}{\partial
x_m}\Big), \qquad \frac{\partial_r f}{\partial \bar
q}=\bar{\partial}_r f  = \frac{1}{2} \Big(\frac{\partial f}{\partial
x_0}+\sum_{m=1}^{3}\frac{\partial f}{\partial x_m}i_m\Big),
\\
\frac{\partial_l f}{\partial q}=\partial_l f & =  \frac{1}{2}
\Big(\frac{\partial f}{\partial x_0}-\sum_{m=1}^{3}i_m\frac{\partial
f}{\partial x_m}\Big), \qquad \frac{\partial_r f}{\partial
q}=\partial_r f  =  \frac{1}{2}\Big(\frac{\partial f}{\partial x_0}
-\sum_{m=1}^{3}\frac{\partial f}{\partial x_m}i_m\Big),
\end{split}\end{equation} $$
\Delta f  =\frac{\partial^2 f}{\partial x_0^2}+\frac{\partial
f^2}{\partial x_1^2}+ \frac{\partial^2 f}{\partial
x_2^2}+\frac{\partial^2 f}{\partial x_3^2}.$$ Note that
$\bar{\partial}_l,\partial_l,\bar{\partial}_r$ and $\partial_r$ all
commute, and
$\Delta=4\partial_r\bar{\partial}_r=4\partial_l\bar{\partial}_l$.
Since the theories generated by $\partial_l$ and $\partial_r$ are
symmetric, we shall use the left operator $\partial_l$ and call the
function simply $q$-holomorphic or regular. If $f:M\to\mathbb H$ a
real-differentiable function, the Cauchy-Riemann-Fueter
equation~\eqref{crf} can be written as $\bar{\partial}_lf=0$ i.~e.
$$\bar{\partial}_l f=\frac{\partial_l f}{\partial \bar
q}=\frac{\partial f}{\partial
(x_0)_k}+\sum_{m=1}^{3}i_m\frac{\partial f}{\partial
(x_m)_k}=0,\quad k=1,\ldots,n.$$

We say that vector fields are {\it tangential Cauchy-Riemann-Fueter
operators} if the following holds.
\begin{itemize}
\item[(i)]{Vector fields on $\partial \mathcal U^1$ arise by restricting (to $\partial \mathcal U^1$) of
vector fields that, in coordinates $(q_1,q_2)$ of $\mathbb H^2$, can
be written in the form
\begin{equation}\label{vf2}\alpha_1\frac{\partial_l}{\partial \bar q_1}+\alpha_2\frac{\partial_l}{\partial
\bar q_2},\end{equation} where the $\alpha_j$ are quaternion valued
functions. This is the class of first-order differential operators
that annihilate $q$-holomorphic functions.}
\item[(ii)]{Vector fields are tangential at $\partial \mathcal U^1$ in the sense that
$$\alpha_1\frac{\partial_l r(q)}{\partial \bar q_1}+\alpha_2\frac{\partial_lr(q)}{\partial \bar q_2}=0,$$
wherever $r(q)=0$. Here $r(q)$ is the defining function given
by~\eqref{r}.}
\end{itemize}

We consider the vector fields on the boundary $\partial\mathcal U_1$
of the Siegel half space with the coordinates
$(q,t)=(x_0,x_1,x_2,x_3,t_1,t_2,t_3)$:
\begin{equation}
\label{vf}
\begin{split}
X_0(q,t)  = &
\partial_{x_{0}}-2x_1\partial_{t_{1}}-2x_2\partial_{t_{2}}-2x_3\partial_{t_{3}},
\\
X_1(q,t) = &
\partial_{x_{1}}+2x_0\partial_{t_{1}}-2x_3\partial_{t_{2}}+2x_2\partial_{t_{3}},
\\
X_2(q,t) = &
\partial_{x_{2}}+2x_3\partial_{t_{1}}+2x_0\partial_{t_{2}}-2x_1\partial_{t_{3}},
\\
X_3(q,t) = &
\partial_{x_{3}}-2x_2\partial_{t_{1}}+2x_1\partial_{t_{2}}+2x_0\partial_{t_{3}}.
\end{split}
\end{equation} They form a basis of the subbundle ${\mathcal T}(\partial\mathcal
U_1)$ of the tangent bundle $T(\partial\mathcal U_1)$. In the same
time the vector fields can be considered as a basis of Lie algebra,
which is an infinitesimal representation of the group $\mathcal Q$.
The commutators between the vector fields~\eqref{vf} is given in
Table~\ref{t2}.
\begin{table}[ht]
\caption{Commutators of the vector fields~\eqref{vf}} \label{t2}
\begin{center}
\begin{tabular}{| c | r | r | r | r |}
\hline \     &  $X_0$ & $X_1$ & $X_2$ & $X_3$
\\  \hline
$X_0$ &  $0$ & $4\partial_{t_1}$ & $4\partial_{t_2}$ &
$4\partial_{t_3}$
\\ \hline
$X_1$ & $-4\partial_{t_1}$ & $0$ & $-4\partial_{t_3}$ &
$4\partial_{t_2}$
\\ \hline
$X_2$ & $-4\partial_{t_2}$ & $4\partial_{t_3}$ & $0$ &
$-4\partial_{t_1}$
\\ \hline
$X_3$ & $-4\partial_{t_3}$ & $-4\partial_{t_2}$ & $4\partial_{t_1}$
& $0$
\\ \hline
\end{tabular}
\end{center}
\end{table}
Let us verify, that the vector fields~\eqref{vf} are tangential
Cauchy-Riemann-Fueter operators. To show this, we  define the
quaternion vector fields
\begin{equation}\label{hh}\begin{split}
\bar H(w,t) & =\frac{1}{2}\big(X_0+i_1X_1+i_2X_2+i_3X_3\big)=\frac{\partial_l}{\partial\bar w}
+wi_1\partial_{t_1}+wi_2\partial_{t_2}+wi_3\partial_{t_3},\\
H(w,t) &
=\frac{1}{2}\big(X_0-i_1X_1-i_2X_2-i_3X_3\big)=\frac{\partial_l}{\partial
w} -i_1\bar w\partial_{t_1}-i_2\bar w\partial_{t_2}-i_3\bar
w\partial_{t_3},
\end{split}\end{equation}
where the terms $wi_k$ and $i_k\bar w$, $k=1,2,3$, are the
quaternion product. Then the commutative relation is $[\bar
H,H]=2\sum_{k=1}^{3}i_k\partial_{t_k}$.

\begin{lemma}\label{lem1}
Using the identification of $\mathcal Q$ with $\partial \mathcal
U_1$, the vector field $\bar H$ is tangential Cauchy-Riemann-Fueter
operator on $\partial \mathcal U_1$.
\end{lemma}
\begin{proof}
We remember the quaternion coordinates $(w,t,r)$ on $\mathcal U_1$:
for $(q_1,q_2)\in \mathcal U_1$
$$w=q_1,\quad t_k=\im_k q_2,\ k=1,2,3,\quad r(q_1,q_2)=\re q_2-|q_1|^2.$$ We use the left multiplication
chain rule $$\frac{\partial_l}{\partial\bar
q_m}=\frac{\partial_l\bar w}{\partial\bar q_m}
\frac{\partial_l}{\partial\bar w}+\frac{\partial_l w}{\partial\bar
q_m}\frac{\partial_l}{\partial w}+\sum_{k=1}^{3}\frac{\partial_l
t_k}{\partial\bar q_m}\frac{\partial}{\partial t}+\frac{\partial_l
r}{\partial\bar q_m}\frac{\partial}{\partial r},\ \ m=1,2.$$
Calculating the derivatives
$$\frac{\partial_l w}{\partial\bar q_m}=\frac{\partial_l \bar w}{\partial\bar q_2}
=\frac{\partial_l t_k}{\partial\bar q_1}=0, \quad\frac{\partial_l
\bar w}{\partial\bar q_1}=1, \quad\frac{\partial_l t_k}{\partial\bar
q_2}=\frac{i_k}{2}, \quad\frac{\partial_l r}{\partial\bar q_1}=-q_1,
\quad \frac{\partial_l r}{\partial\bar q_2}=\frac{1}{2},$$ for
$k=1,2,3$, $m=1,2$, we conclude
$$\frac{\partial_l}{\partial\bar
q_1}=\frac{\partial_l}{\partial\bar w}-q_1
\frac{\partial_l}{\partial r},\qquad \frac{\partial_l}{\partial\bar
q_2}=\sum_{k=1}^{3}\frac{i_k}{2} \frac{\partial_l}{\partial\bar
t_k}+\frac{1}{2}\frac{\partial_l}{\partial r}. $$ Multiplying the
second equation by $2q_1$ from the left and summing the derivative,
we get
$$\frac{\partial_l}{\partial\bar w}+\sum_{m=1}^{3}wi_m\partial_{t_m}
=\frac{\partial_l}{\partial\bar
q_1}+2q_1\frac{\partial_l}{\partial\bar q_2}.$$ Notice that
$\frac{\partial_l r}{\partial\bar q_1}+2q_1\frac{\partial_l
r}{\partial\bar q_2}=0$. This proves that the operator
$\frac{\partial_l}{\partial \bar
q_1}+2q_1\frac{\partial_l}{\partial\bar q_2}$ is tangential
Cauchy-Riemann-Fueter operator.
\end{proof}

\section{The Cauchy-Szego kernel}

We will need the following $3$-form $Dq$, which is different from
the usual volume $3$-form in~$\mathbb R^3$. We write $v=dx_0\wedge
dx_1\wedge dx_2\wedge dx_3$ for the volume form in $\mathbb H$. The
form $Dq$ is defined as an alternating $\mathbb R$ trilinear
function by \begin{equation} \label{Dform} \langle
h_1,Dq(h_2,h_3,h_4)\rangle =v(h_1,h_2,h_3,h_4)
\end{equation} for all $h_1,\ldots,h_4\in\mathbb H$. The coordinate expression for $Dq$ is
$$Dq=dx_1\wedge dx_2\wedge dx_3-i_1dx_0\wedge dx_2\wedge dx_3-i_2dx_0\wedge dx_3\wedge dx_1
-i_3dx_0\wedge dx_1\wedge dx_2.$$ Other properties of $Dq$ reader
can fined in~\cite{Sud}. This is the principal form for the
formulations of the integral theorems related to the $q$-holomorphic
functions. We present the Cauchy integral theorem here.

\begin{theorem}[\cite{Sud}]\label{cth}
Suppose $f$ is $q$-holomorphic function in an open set $U\in\mathbb
H$. Let $q_0$ be a point in $U$, and let $C$ be a rectifiable
$3$-chain which is homologous, in the singular homology of
$U\setminus \{q_0\}$, to a differentiable $3$-chain whose image is
$\partial B$ for some ball $B\subset U$. Then
$$\frac{1}{2\pi^2}\int_{C}\frac{(q-q_0)^{-1}}{|q-q_0|^2}Dqf(q)=nf(q_0),$$
where $n$ is the wrapping number of $C$ about $q_0$.
\end{theorem}

We start from the definition of the Hardy space $\mathcal
H^2(\mathcal U_1)$. In the sequel, we write $\ii\cdot
t=\sum_{m=1}^{3}i_mt_m$. Let $dh$ be Haar measure on $\mathcal Q$.
Using the identification of $\partial \mathcal U_1$ with $\mathcal
Q$ we introduce the measure $d\beta$ on $\partial \mathcal U_1$. So,
the integration formula $$\int_{\partial \mathcal
U_1}F(q)d\beta(q)=\int_{\mathbb H\times \mathbb R^3}
F(q_1,|q_1|^2+\ii\cdot t)\,dq_1\,dt,$$ holds for a continuous
function $F$ of compact support. With this measure we can define the
space $L^2(\mathcal Q)=L^2(\partial \mathcal U_1)$.

For any function $F$ defined on $\mathcal U_1$, we write
$F_{\varepsilon}$ for its "vertical translate" (we mean that the
vertical direction is given by the positive direction of $\re q_2$):  $$F_{\varepsilon}(q)=F(q+\varepsilon \mathbf e),
\quad\text{where}\quad \mathbf e=(0,0,0,0,1,0,0,0).$$ If
$\varepsilon>0$, then $F_{\varepsilon}$ is defined in the
neighborhood of $\partial \mathcal U_1$. In particular,
$F_{\varepsilon}$ is defined on $\partial \mathcal U_1$.

\begin{definition} The space $\mathcal H^2( \mathcal U_1)$ consists of all functions $F$
holomorphic on $\mathcal U_1$, for which
\begin{equation}\label{hs}\sup_{\varepsilon>0}\int_{\partial
\mathcal U_1}|F_{\varepsilon}(q)|^2\,
d\beta(q)<\infty.\end{equation} \end{definition} The norm
$\|F\|_{\mathcal H^2( \mathcal U_1)}$ of $F$ is then the square root
of the left-hand side of~\eqref{hs}. The space $\mathcal H^2(
\mathcal U_1)$ is a Hilbert space.

We shall need the following statements concerning the holomorphic
$\mathcal H^2$ space on the upper half space $\mathcal U_1$ of
 $\mathbb H^2$. We denote the upper half space $\mathcal U_1$ by $\mathbb R^{4}_+$.
 Define $\mathcal H^2(\mathbb R^{4}_+)$
 to be a set of functions $f(q)$ which are $q$-holomorphic for $q=u+\ii\cdot v$ in the upper half space $u>0$
satisfying \begin{equation} \label{eq1} \Big(\sup_{u>0}\int_{\mathbb
R^{3}}|f(u+\ii\cdot v)|^2\,dv\Big)^{1/2}=\|f\|_{\mathcal H^2(\mathbb
R^{4}_+)}<\infty.
\end{equation}  The {\it maximal inequality} \begin{equation}
\label{eq2} \int_{\mathbb R^{3}}\sup_{u>0}|f(u+\ii\cdot
v)|^2\,dv\leq A^2\|f\|^2_{\mathcal H^2(\mathbb R^{4}_+)}
\end{equation} is a consequence of~\eqref{eq1}. This implies that there exists an $f^{b}\in L^2(\mathbb R^3)$
such that \begin{equation} \label{eq4}f(u+\ii\cdot v)\to
f^b(v)\quad\text{as}\quad u\to 0,\quad\text{for a.~e.}\quad
v,\end{equation} and the identity \begin{equation} \label{eq3}
\|f^{b}\|_{L^2(\mathbb R^3)}=\|f\|_{\mathcal H^2(\mathbb R^{4}_+)}.
\end{equation}
The definition of $\mathcal H^2({\mathbb R}^{4}_+)$ and further
properties can be found in~\cite{St}.

Let $q_1\in\mathbb H$ and $\delta>0$. We define $f$ on $\mathbb
R^{4}_+$ by the following $f(q_2)=F(q_1,q_2+\delta +|q_1|^2)$.
\begin{lemma}\label{1bv}
Suppose that $F\in\mathcal H^2(\mathcal U_1)$ and $\delta>0$. Then
for each $q_1\in\mathbb H$, the function $$f(q_2)=F(q_1,q_2+\delta
+|q_1|^2)$$ with $q_2=u+\ii\cdot v$, belongs to $\mathcal
H^2({\mathbb R}^{4}_+)$
\end{lemma}
\begin{proof}
We use the action~\eqref{acgr} of group $\mathcal Q$ on $\mathcal
U_1$ and others $q$-holomorphic automorphisms of $\mathcal U_1$ to
simplify the statement. Let $h=[w,t]\in\mathcal Q$ and
$q=(q_1,q_2)$. It is easy to see that $h(q+\varepsilon \mathbf
e)=h(q)+\varepsilon \mathbf e$. Then $F(h(q))\in\mathcal
H^2(\mathcal U_1)$ whenever $F(q)\in\mathcal H^2(\mathcal U_1)$,
since the measure $d\beta(q)$ is invariant with respect to action of
$h$. We take $h=[q_1,0]$ and reduce the considerations to the
function $f(q_2)=F(0,q_2+\delta)$. Applying the dilations and the
property $d\beta(\delta(q))=\delta^{10}d\beta$, we can consider
$\delta=1$ and work with $$f(u+\ii\cdot v)=F(0,u+1+\ii\cdot
v),\quad\text{for}\quad u>0.$$ By the standard way, from the Cauchy
integral Theorem~\ref{cth} one can obtain the mean value property
for $q$-holomorphic functions \begin{equation} \label{mvp}
|f(u+\ii\cdot v)|^2=|F(0,u+1+\ii\cdot v)|^2\leq
c\int\limits_{|q_1|^2+|q_2|^2< 1/4}|F(q_1,q_2+u+1+\ii\cdot v)|^2\,
dq_1\,dq_2.
\end{equation} where $c^{-1}=\frac{\pi^4}{2^84!}$ is the volume of the ball of the radius $1/2$ in $\mathbb H^2$.
Notice that since $|\re q_2|<1/2$ and $u>0$, we have
$\re(q_2+u+1)\in(1/2,3/2)$. Thus $\re(q_2+u+1)>|q_1|^2$ and this
guarantees that the integration in~\eqref{mvp} is taken over some
subset of  $\mathcal U_1$. We write $q_2=x+\ii\cdot y$ and
integrate~\eqref{mvp} with respect to $v$ over $\mathbb R^3$. After
applying the Fubini theorem, we deduce $$\int_{\mathbb
R^3}|f(u+\ii\cdot v)|^2\,dv\leq c\int\limits_{|q_1|\leq 1/2}
\int\limits_{\mathbb
R^3}\int\limits_{|q_2|<1/2}|F(q_1,q_2+u+1+\ii\cdot
v)|^2\,Dq_2\,dx\,dq_1\,dv.$$ We take now the second integral with
respect to $Dq_2$ and write $v=y$. We get $$\int_{\mathbb
R^3}|f(u+\ii\cdot v)|^2\,dv\leq c_1\int\limits_{|q_1|<
1/2}\int\limits_{\mathbb
R^3}\int\limits_{|x|<1/2}|F(q_1,x+u+1+\ii\cdot
y)|^2\,dx\,dq_1\,dy.$$ We make the following change of variables
$x+u+1=\varepsilon +|q_1|^2$. Since $|x|<1/2$ and $|q_1|^2<1/4$, the
rang of new variable $\varepsilon$ is in the interval
$(u+1/4,u+3/2)$. So, the last integral take the form
\begin{equation*}\begin{split}&
\int\limits_{u+3/2}^{u+1/4}\int\limits_{|q_1|<
1/2}\int\limits_{\mathbb R^3}|F(q_1,\varepsilon +|q_1|^2+\ii\cdot
y)|^2\,dq_1\,dy\,d\varepsilon
\\ & \leq\frac{5}{4}\int\limits_{\mathbb H\times\mathbb R^3}|F(q_1,\varepsilon +|q_1|^2+\ii\cdot y)|^2\,dq_1\,dy
=\frac{5}{4}\int\limits_{\partial\mathcal
U_1}|F(q+\varepsilon\mathbf e)|^2\,d\beta(q)<\infty
\end{split}\end{equation*} because of~\eqref{hs}. This shows that $f\in\mathcal H^2(\mathbb R^{4}_+)$
and we finish the proof of Lemma~\ref{1bv}.
\end{proof}

\begin{theorem}
\label{bv} Suppose $F$ belongs to $\mathcal H^2( \mathcal U_1)$.
Then \begin{itemize}
  \item[1.]{There exists an $F^b\in L^2(\partial \mathcal U_1)$ so that $F(z+\varepsilon \mathbf e)
  \vert_{\partial \mathcal U_1}\to F^b$ in the $L^2(\partial \mathcal U_1)$ norm, as $\varepsilon\to 0$.}
  \item[2.]{The above mentioned space of $F^b$ is a closed subspace of $L^2(\partial \mathcal U_1)$.
  Moreover }
  \item[3.]{$\|F^b\|_{L^2(\partial\mathcal U_1)}=\|F\|_{\mathcal H^2( \mathcal U_1)}$.}
\end{itemize}
\end{theorem}

\begin{proof} To proof Theorem~\ref{bv} we exploit the reducing to the one dimensional case.
Let us fix $q_1$ and consider the function $F(q)=f(q_1,q_2)$ as a
function of $q_2$. We apply the maximal inequality~\eqref{eq2} and
equality~\eqref{eq3} to the function $f(q_2)=F(q_1,
q_2+\delta+|q_1|^2)$. We continue to write $q_2=x+\ii\cdot y$. Then
\begin{eqnarray*} \int_{\mathbb R^3}\sup_{x>0}|F(q_1,
x+\delta+|q_1|^2+\ii\cdot y)|^2\,dy\leq A^2\int_{\mathbb
R^3}|F(q_1,\delta+|q_1|^2 +\ii\cdot y)|^2\,dy.
\end{eqnarray*} We integrate over $q_1\in\mathbb H$, write $x=\varepsilon$, use the definition of $d\beta=dq_1\,dy$,
and $|q_1|^2=\re q_2$ on the $\partial\mathcal U_1$. Then the last
inequality can be written as $$\int_{\partial\mathcal
U_1}\sup_{\varepsilon>0}|F(q+
\varepsilon+\delta)|^2\,d\beta(q)\leq A^2\int_{\partial\mathcal
U_1}|F(q+\delta)|^2\,d\beta(q).$$ Letting $\delta$ to $0$, we see
$$\int_{\partial\mathcal U_1}\sup_{\varepsilon>0}|F(q+
\varepsilon)|^2\,d\beta(q)\leq A^2\|F\|^2_{\mathcal H^2(\mathcal
U_1)}.$$ We conclude from the last inequality that for almost all
$q_1\in\mathbb H$ the function $F(q_1,q_2+|q_1|^2)$, as a function
of $q_2$, is in $\mathcal H^2({\mathbb R}^{4}_+)$ and
by~\eqref{eq4} we see that the limit $\lim_{\varepsilon\to 0}
F(q+\varepsilon\mathbf e)=F^b(q)$ exists for almost every
$q\in\partial\mathcal U_1$.

We show now the property~3. By Fatou's lemma, it follows
$$\int_{\partial\mathcal U_1}|F^b(q)|^2\,d\beta(q)\leq
\sup_{\varepsilon>0}\int_{\partial\mathcal U_1}|F(q+
\varepsilon\mathbf e)|^2\,d\beta(q) =\|F\|^2_{\mathcal H^2(\mathcal
U_1)}.$$ From the other hand, we use~\eqref{eq3} and get
$$\int_{\mathbb R^3}|F(q_1,\varepsilon+|q_1|^2+\ii\cdot
y))|^2\,dy\leq \sup_{\varepsilon>0} \int_{\mathbb
R^3}|F(q_1,\varepsilon+|q_1|^2+\ii\cdot y))|^2\,dy=\int_{\mathbb
R^3}|F(q_1,|q_1|^2+\ii\cdot y))|^2\,dy.$$ Integrating over
$q_1\in\mathbb H$ both sides of the last inequality and taking
supremum over $\varepsilon >0$, we get the property~3.

We showed that the limit $F^b$ in $\mathcal H^2(\mathcal U_1)$ norm
exists. To complete the proof of the assertion ~2 we note that the
mean value property~\eqref{mvp} says that for any compact
$K\in\mathcal U_1$ there is a constant $c_K$, such that
\begin{equation}\label{comp}\sup_{q\in K}|F(q)|\leq
c_K\|F\|_{\mathcal H^2(\mathcal U_1)}.\end{equation} We conclude,
that if a sequence $F_n$ converges in $\mathcal H^2(\mathcal U_1)$
norm, then the sequence $F_n$ converges uniformly on compact subsets
of $\mathcal U_1$, that implies that the space $\mathcal
H^2(\mathcal U_1)$ is complete with respect of its norm. The
Theorem~\ref{bv} is proved.\end{proof}

We determine now the Cauchy-Szeg\"o kernel $S(q,\omega)$ for the
domain $\mathcal U_1$. The Cauchy-Szeg\"o kernel $S(q,\omega)$ is a
quaternion valued function, defined on $\mathcal U_1\times\mathcal
U_1$ and satisfying the following conditions
\begin{itemize}
  \item[1.]{For each $\omega\in \mathcal U_1$, the function $q\mapsto S(q,\omega)$ is $q$-holomorphic for $q\in
  \mathcal U_1$, and belongs to $\mathcal H^2(\mathcal U_1)$. This allows us to define, for each $\omega\in
  \mathcal U_1$, the boundary value $S^b(q,\omega)$
  for almost all
$q\in\partial \mathcal U_1$.}
  \item[2.]{The kernel $S$ is symmetric: $S(q,\omega)=\overline{S(\omega,q)}$ for each $(q,\omega)\in
  \mathcal U_1\times\mathcal U_1$. The symmetry permit us to extend the definition of  $S(\omega,q)$
  so that for each $q\in\mathcal U_1$,  the function $S^b(q,\omega)$ is defined for almost every
  $\omega \in \partial \mathcal U_1$.}
  \item[3.]{The kernel $S$ satisfies the reproducing property in the following sense \begin{equation}
\label{rp} F(q)=\int_{\partial \mathcal
U_1}S(q,\omega)F^b(\omega)\,d\beta(\omega),\qquad q\in\mathcal U_1,
\end{equation} whenever $F\in\mathcal H^2(\mathcal U_1)$.}
\end{itemize}
\begin{theorem}\label{csk}
Let $S(q,\omega)=\frac{k}{r^5(q,\omega)}$, where \begin{equation}
\label{ r1} r(q,\omega)=\frac{q_2+\bar\omega_2}{2}-\bar\omega_1
q_1\quad\text{and}\quad k=\frac{3}{8\pi^{4}}.
\end{equation} Then $S(q,\omega)$ is the (unique) function that satisfies the properties~$1.-3.$ above.
\end{theorem}
\begin{proof}
We introduce the Cauchy-Szeg\"o projection operator $C$. The
operator $C$ is the orthogonal projection from $L^2(\partial\mathcal
U_1)$ to the subspace of functions $\{F^b\}$ that are boundary
values of functions $F\in\mathcal H^2(\mathcal U_1)$. So for each
$f\in L^2(\partial\mathcal U_1)$, we have that $C(f)=F^b$ for some
$F\in\mathcal H^2(\mathcal U_1)$; moreover, $C(F^b)=F^b$ and $C$ is
self-adjoint: $C^*=C$.

We fix $q\in \mathcal U_1$. Then $(Cf)(q)=F(q)$, where $F$
corresponds to $F^b$. The kernel $S(q,\omega)$ will be defined by
the representation \begin{equation} \label{csrep}
F(q)=\int_{\partial \mathcal
U_1}S(q,\omega)f(\omega)\,d\beta(\omega).
\end{equation}

The space $\mathcal H^2(\mathcal U_1)$ is identified with some
subspace of $L^2(\partial \mathcal U_1)$. We take a basis
$\{\varphi_j\}$ of this subspace. We can expect that
\begin{equation}\label{bas}S(q,\omega)
=\sum_{j}\varphi_j(q)\bar\varphi_j(\omega).\end{equation} (In this
case we also have the symmetry property
$$\overline{S(q,\omega)}=\sum_{j}\overline{\varphi_j(q)\bar\varphi_j(\omega)}
=\sum_{j}\varphi_j(\omega)\bar\varphi_j(q)=S(\omega,q).\Big)$$
Indeed, if $\sum_{j}|a_j|^2<\infty$, then
$\sum_ja_j\varphi_j\in\mathcal H^2(\mathcal U_1)$ and
$\|\sum_ja_j\varphi_j\|_{\mathcal H^2(\mathcal U_1)}
=\sum_{j}|a_j|^2$. Moreover, for any compact set $K\subset \mathcal
U_1$ we have by~\eqref{comp} that
$$\sup\limits_{q\in K}\Big|\sum_ja_j\varphi_j(q)\Big|\leq
c_K\Big(\sum_{j}|a_j|^2\Big)^{1/2}.$$ Applying the converse of
Schwart's inequality, we get
$$\Big(\sum_{j}|\varphi_j(q)|^2\Big)^{1/2}\leq c_K, \quad\text{for
all}\quad q\in K.$$ Thus the sum~\eqref{bas} converges uniformly
whenever $(q,\omega)$ belongs to a compact subset of $\mathcal
U_1\times \mathcal U_1$.

Let us take~\eqref{bas} as the definition of $S(q,\omega)$. By the
symmetry $\overline{S(q,\omega)}=S(\omega,q)$ the function
$\overline{S(q,\omega)}$ belongs to $\mathcal H^2(\mathcal U_1)$ for
each fixed $q\in \mathcal U_1$. Moreover $\overline{S(q,\omega)}$
extends to $q\in \mathcal U_1$, $\omega\in\partial \mathcal U_1$, by
the identity~\eqref{bas} with the series converging in the norm of
$L^2(\partial \mathcal U_1)$.

These arguments establish the existence of the function $S$
satisfying~\eqref{csrep} and the properties~1.-3. The reproducing
property~\eqref{rp} uniquely determines $\overline{S(q,\omega)}$ as
an element of $\mathcal H^2(\mathcal U_1)$ for each fixed $q$.
Together with conclusion~3 of Theorem~\ref{bv}, this shows that $S$
is uniquely determined by the properties 1.-3.

We continue the proof establishing the precise shape of  the
function $S(q,\omega)$. We use the translation $h(q)$, the unitary
rotation $\mathcal R(q)$ and the dilation $\delta(q)$ on $\mathcal
U_1$. Observe that the measure $d\beta$  is invariant with respect
to translation and unitary rotations and
$d\beta(\delta(q))=\delta^{10}d\beta$, where $10$ is the homogeneous
dimension of the group $\mathcal Q$. Notice that the space $\mathcal
H^2(\mathcal U_1)$ is also takes into itself under the above
mentioned transformations. Then, we have
$$F(q)=\int_{\partial \mathcal U_1}S(h(q),h(\omega))F^b(\omega)\,d\beta(\omega),$$
$$F(q)
=\int_{\partial \mathcal U_1}S(\mathcal R(q),\mathcal
R(\omega))F^b(\omega)\,d\beta(\omega),$$
$$F(q)=\int_{\partial \mathcal U_1}S(\delta(q),\delta(\omega))\delta^{10}F^b(\omega)\,d\beta(\omega).$$
We conclude
\begin{equation}\label{sym}S(q,\omega)=S(h(q),h(\omega)),\quad
S(q,\omega) =S(\mathcal R(q),\mathcal R(\omega)),\quad
S(q,\omega)=S(\delta(q),\delta(\omega))\delta^{10}.\end{equation}
The identities~\eqref{sym} hold for each $(q,\omega)\in \mathcal
U_1\times\mathcal U_1$. They also hold for each $q\in\mathcal U_1$
and almost all $\omega\in\partial \mathcal U_1$, because for each
fixed $q\in\mathcal U_1$ the identities~\eqref{sym} are the
identities for the elements of $L^2(\partial \mathcal U_1)$ and,
therefore, they are fulfilled for almost all $\omega\in\partial
\mathcal U_1$.

Let $S(q)=S(q,0)$, then $S(q)$ is holomorphic on $\mathcal U_1$ and
is independent of $q_1$. So, we may write $S(q)=s(q_2)$. Using the
transformation of $S$ with respect to dilation we get
$s(\delta^2q_2)=\delta^{-10}s(q_2)$ for all positive $\delta$. We
conclude that $S(q)=cq_2^{-5}$. We use now the translation~$h(q)$.
Let $q\in\mathcal U_1$ and $\omega\in \partial\mathcal U_1$. Then by
identification of $\mathcal U_1$ with $\mathcal Q$ we have
$\omega=h(0)$ and $S(q,h(0))=S(h^{-1}(q),0)=S(h^{-1}(q))$. We
calculate the element $h$ such that $h(0)=(\omega_1,\omega_2)$. We
have
$$h=[w,t]=[\omega_1,\im_1(\omega_2-|\omega_1|^2),\im_2(\omega_2-|\omega_1|^2),\im_3(\omega_2-|\omega_1|^2)].$$
We also observe, that $\re\omega_2=|\omega_1|^2$. Then the second
coordinate of the action $h^{-1}$ is
$$q_2+|\omega_1|^2-\ii\cdot\im
\omega_2-2\bar\omega_1q_1=q_2+\bar
\omega_2-2\bar\omega_1q_1=2r(q,\omega)$$ by~\eqref{acgr}. Therefore
$S(q,\omega)=kr^{-5}(q,\omega)$ with $k=2^{-5}c$.

We calculate the constant $k$. Use the reproducing formula for the
function $F(q)=(q_2+1)^{-5}=(2\overline{r(\mathbf e,q)})^{-5}$.
Applying the reproducing formula, we get \begin{eqnarray*} 2^{-5} &
= & F(\mathbf e)=k\int_{\partial\mathcal U_1}
r(e,q)^{-5}F(q)\,d\beta{q}
\\
& =  &  2^{5} k\int_{\partial\mathcal U_1} |F(q)|^2\,d\beta{q}=2^{5}
k\int_{q^{\prime}\in\mathbb H} \int_{t\in\mathbb R^3}|F(q^{\prime},
|q^{\prime}|^2+\sum_{m=1}^{3}i_mt_m)|^2\,dq^{\prime}dt.
\end{eqnarray*} Thus \begin{equation}\label{in}k^{-1}=4^{5}\int_{w\in\mathbb H}\,dw
\int_{t\in\mathbb R^3}\big(|t|^2+(|w|^2+1)^2\big)^{-5}dt.
\end{equation}
Taking the integral over $\mathbb R^3$, we observe that~\eqref{in}
is reduced to the product of four values $\alpha,\beta,\gamma,
\delta$, where $\alpha=4\pi$ is the volume of the unite sphere in
$\mathbb R^3$,
$$\beta=2\pi^{2}\quad\text{is the volume of
the unite sphere in}\quad \mathbb H,$$
$$\gamma=\int_{0}^{\infty}\frac{r^2\,dr}{(r^2+1)^{5}}=\frac{\Gamma(3/2)\Gamma(7/2)}{2\Gamma(5)},
\quad\text{and}\quad$$
$$\delta=\int_{0}^{\infty}\frac{\rho^{3}\,d\rho}{(\rho^2+1)^{7}}=\frac{\Gamma(5)}{2\Gamma(7)}.$$
Multiplying all values, we get $k=\frac{3}{8\pi^{4}}$. \end{proof}

Notice the following. \begin{itemize}
  \item[(i)]{The function $r(q,\omega)$ is $q$-holomorphic in $q$ and anti $q$-holomorphic in $\omega$.
  If $q=\omega$, the function $r$ agrees with the function~\eqref{r}.}
  \item[(ii)]{For each fixed $\omega\in\mathcal U_1$ the function $r(q,\omega)$ (and hence $S(q,\omega)$)
  is $q$-holomorphic for $q$ in a neighborhood of closure of $\mathcal U_1$. For each fixed $q\in\mathcal U_1$
  the function $r(q,\omega)$ (and hence $S(q,\omega)$) is anti $q$-holomorphic for $\omega$ in a neighborhood of
  closure of $\mathcal U_1$. In particular, if $q\in \mathcal U_1$ is fixed, the boundary function
  $S^b(\cdot,\omega)$ is, actually, defined on all of the boundary~$\partial\mathcal U_1$.}
\end{itemize}

\section{The projection operator}

We describe the projection operator $C$ as a convolution operator on
the group $\mathcal Q$. We remember that mapping $f\mapsto C(f)$
assigns to each element $f\in L^2(\partial\mathcal U_1)$ another
element of $L^2(\partial\mathcal U_1)$ of the form $C(f)=F^b$, for
some $F\in\mathcal H^2(\mathcal U_1)$. Theorem~\ref{bv} and
reproducing property~\eqref{csrep} imply
\begin{equation}\label{cf}(Cf)(q)=\lim_{\varepsilon\to
0}F(q+\varepsilon\mathbb e)\vert_{\partial\mathcal
U_1}=\lim_{\varepsilon\to 0}\int_{\partial\mathcal
U_1}S(q+\varepsilon\mathbb
e,\omega)f(\omega)\,d\beta(\omega),\end{equation} where $q\in
\partial\mathcal U_1$ and the limit is taken in
$L^2(\partial\mathcal U_1)$ norm.

We now use the identification of $\partial\mathcal U_1$ with
$\mathcal Q$. We write $\partial\mathcal U_1\ni\omega=g(0)$ for the
unique $g\in\mathcal Q$, $\partial\mathcal U_1\ni q=h(0)$ for
$h\in\mathcal Q$, and $d\beta(\omega)=dg$. The
properties~\eqref{sym} of the Cauchy-Szeg\"o kernel give
$$S(q+\varepsilon\mathbb e,\omega)=S(g^{-1}(q+\varepsilon\mathbb
e),g^{-1}(\omega))=S(g^{-1}(q)+\varepsilon\mathbb
e,g^{-1}(\omega))=S(g^{-1}(h(0))+\varepsilon\mathbb e,0).$$ We
change the notation setting
$K_{\varepsilon}(h)=S(g^{-1}(h(0))+\varepsilon\mathbb e,0)$,
$f(g)=f(g(0))=f(\omega)$, and $(Cf)(h)=(Cf)(h(0))=(Cf)(q)$.
Then~\eqref{cf} takes the form \begin{equation} \label{cf1}
(Cf)(h)=\lim_{\varepsilon\to 0}\int_{\mathcal
Q}K_{\varepsilon}(g^{-1}\circ h)f(g)\,d\beta(g),
\end{equation} for $f\in L^2(\mathcal Q)$, where the limit is taken in $L^2(\mathcal Q)$.
We see that~\eqref{cf1} is formally written as the convolution $(Cf)(h)=(K\ast f)(h)$,
where $K=\lim_{\varepsilon\to 0}K_{\varepsilon}$ is a distribution.

We obtain the precise  form of $K(h)$. Let $h=[w,t]$, then by
identification of $\partial\mathcal U_1$ with $\mathcal Q$ and
Theorem~\ref{csk} we have
$K_{\varepsilon}(h)=c(|w|^2+\varepsilon+\ii\cdot t)^{-5}$ with
$c=\frac{6}{\pi^{4}}$. We observe that
$$K_{\varepsilon}(h)=\frac{-c}{2i_13i_24i_3}\frac{\partial^3}{\partial
t_3\partial t_2\partial t_1}\Big((|w|^2+\varepsilon+\ii\cdot
t)^{-2}\Big).$$ The function $(|w|^2+\varepsilon+\ii\cdot t)^{-2}$
is locally integrable on $\mathcal Q$ (see remark below). Passing to
the limit, we see that the distribution $K$ is given by
$$K(h)=\frac{c}{24}\frac{\partial^3}{\partial t_3\partial
t_2\partial t_1}\Big((|w|^2+\ii\cdot t)^{-2}\Big)$$ and equals to
the function $(|w|^2+\ii\cdot t)^{-5}$ away from the origin.

\begin{remark}
On the group $\mathcal Q$ we can introduce the homogeneous norm
$$\|h\|^2=\|[w,t]\|^2=|w|^2+|t|,$$ where $|\cdot|$ denotes the
Euclidean norm. The homogeneous norm is a homogeneous of order $1$
with respect to dilation $\delta$  function, namely: $\|\delta
(h)\|=\delta\|h\|$. We also recall the analogue of the integration
formula in the polar coordinates~\cite{FS}. There is a positive
constant $k$, such that whenever $f$ is a non-negative function on
$(0,\infty)$, then
$$\int\limits_{\mathcal Q}f(\|h\|)\,dh=k\int_{0}^{\infty}f(r)r^{Q-1}\,dr,$$ where $Q$ is the homogeneous
dimension of the group $\mathcal Q$.
The homogeneous dimension of $\mathcal Q$ equals $10$. Since the
kernel of the Cauchy-Szeg\"o projection satisfies $|K(h)|\approx
\|h\|^{-10}$ then the function $(|w|^2+\ii\cdot t)^{-2}$ is locally
integrable.  \end{remark}

\begin{remark}
We can change the arguments that we used for the calculation of the precise
form of the distribution $K(h)$. It is sufficient to
differentiate only one time to obtain the integrable function. We
can consider
$$K_{\varepsilon}(h)=\frac{-c}{4i_m}\frac{\partial}{\partial
t_m}\Big((|w|^2+\varepsilon+\ii\cdot t)^{-4}\Big)\quad\text{for
any}\quad m=1,2,3$$ and then argue as above.\end{remark}

\section{The kernel of the operator $\Delta_{\lambda}$}

Recall from Section~3, the operator $$\bar
H(w,t)=\frac{1}{2}\big(X_0+i_1X_1+i_2X_2+i_3X_3\big)
=\frac{\partial_l}{\partial\bar w}
+wi_1\partial_{t_1}+wi_2\partial_{t_2}+wi_3\partial_{t_3},$$ where
the terms $wi_k$, $k=1,2,3$, are the quaternion product. The
conjugate operator is
$$H(w,t)=\frac{1}{2}\big(X_0-i_1X_1-i_2X_2-i_3X_3\big)=\frac{\partial_l}{\partial w}
-i_1\bar w\partial_{t_1}-i_2\bar w\partial_{t_2}-i_3\bar
w\partial_{t_3}.$$ Then the commutative relation give us $[\bar
H,H]=-2\sum_{k=1}^{3}i_k\partial_{t_k}$.

To define $\bar \partial_b$ operator on $\mathcal Q$ we take a
function $f$ and set
$$\bar \partial_b f=\bar Hf\,d\bar q,$$ where $d\bar q=dx_0-\sum_{k=1}^{3}i_kdx_k$.
It is easy to see that $$d\bar q\wedge d\bar q=i_1dx_2\wedge
dx_3+i_2dx_3\wedge dx_1+ i_3dx_1\wedge dx_2.$$
Then
$$\bar
\partial_b^2 f=\bar H\bar H fd\bar q\wedge d\bar q\neq 0.$$
This is the difference with the Heisenberg group. We then define
the formal adjoint $\bar \partial_b^{*}$ of $\bar \partial_b$ by
$$
\bar \partial_b^{*}f=-Hf\quad\text{where}\quad f=f\,d\bar q\quad\text{is a $(0,1)$-form}.
$$
Finally, we can define $\square_b$ operator by
$$\square_b=\bar \partial_b\, \bar \partial_b^*+\bar \partial_b^* \bar \partial_b.$$
When $f$ is a function, then
\begin{equation*}
\begin{split}
\bar\partial_b^*\bar\partial_b (f)  & =-H\bar H(f),\\
\bar\partial_b\bar\partial_b^* (f)  & = 0.
\end{split}
\end{equation*}
It follows that
\begin{equation*}
\begin{split}
-H\bar H&=-\frac{1}{2}(H\bar H+\bar HH)+\frac{1}{2}[\bar H,H]
=-\frac{1}{2}(H\bar H+\bar HH)-2\sum_{k=1}^{3}i_k\partial_{t_k}\\
&=-\frac{1}{4}\Big(\sum_{l=0}^3X_l^2+8\sum_{k=1}^{3}i_k\partial_{t_k}\Big),
\end{split}
\end{equation*}
where $X_l$, $l=0,1,2,3$ are defined in~\eqref{vf} and
$\partial_{t_k}$, $k=1,2,3$ are their commutators. In fact, we may
define a little bit more general operator $\Delta_\lambda$ as
follows
$$\Delta_{\lambda}=\sum_{l=0}^{3}X_l^2+4\sum_{k=1}^{3}i_k\lambda_k\partial_{t_k}.$$
The operator $\Delta_\lambda$ is called the subLaplacian in the
literature and from a theorem of H\"ormander ~\cite{Hor}, it is a
subelliptic operator. In terms of coordinates on ${\mathcal Q}$,
one has
$$\Delta_{\lambda}=\sum_{l=0}^{3}\partial_{x_l}^2+4|x|^2\sum_{k=1}^{3}\partial_{t_k}^2
+4\sum_{k=1}^{3}\big((xi_k \cdot\partial_x)+\lambda_k
i_k\big)\partial_{t_k}.$$
Let us observe that $\Delta_{\lambda}$
possesses the following symmetry properties
\begin{itemize}\item[(i)]{is left invariant on $\mathcal Q$ with  respect to the translation
defined by the group multiplication from the left,}\item[(ii)]{has
degree $2$ with respect to the dilation $\delta$,}\item[(iii)]{is
invariant under the unitary rotation $\mathcal R$ on ${\mathbb
R}^4$.}
\end{itemize}

Recall that partial Fourier transform in the variables $t_k$,
$k=1,2,3$, is defined as follows $$\mathcal
F(f)(x,\tau_1,\tau_2,\tau_3)=\widetilde f
(x,\tau_1,\tau_2,\tau_3)=\int_{\mathbb
R^3}e^{-\sum_{k=1}^{3}i_k\tau_k t_k}f(t)\,dt,$$
$$f(x,t_1,t_2,t_3)=\frac{1}{(2\pi)^3}\int_{\mathbb
R^3}e^{\sum_{k=1}^{3}i_kt_k\tau_k}f(\tau)\,d\tau.$$ Then
$$\mathcal F\Big(\frac{\partial f}{\partial
t_k}\Big)(x,\tau_1,\tau_2,\tau_3)=i_k\tau_k\widetilde
f(\tau),\qquad\mathcal F\Big(\frac{\partial^2 f}{\partial
t_k^2}\Big)(x,\tau_1,\tau_2,\tau_3)=-\tau_k^2\widetilde f(\tau).$$
Consequently,
$$\widetilde{\Delta_{\lambda}}(\tau)=\sum_{l=1}^{4}\partial_{x_l}^2-4|x|^2\sum_{k=1}^{3}\tau_k^2
+4\sum_{k=1}^{3}\big((xi_k
\cdot\partial_x)i_k-\lambda_k\big)\tau_k.$$ We write
$|\tau|^2=\sum_{k=1}^{3}\tau_k^2$. When $\lambda=0$, one has
$$\widetilde{\Delta_{0}}(\tau)=\sum_{l=1}^{4}\partial_{x_l}^2-4|x|^2|\tau|^2
+4\sum_{k=1}^{3}\big(i_k( xi_k\cdot\partial_x)\big)\partial_{t_k}.
$$

The property~(iii) implies that
$\widetilde{\Delta_{\lambda}}(\tau)$ is invariant under the
unitary rotation group action on $\R^4$. Therefore, the
fundamental solution $\widetilde K_\lambda(x,\tau)$ of
$\widetilde{\Delta_{\lambda}}(\tau)$ is a radial distribution. It
follows that
$$
\sum_{k=1}^3\big(xi_k \cdot\partial_x\big)i_k\tau_k\widetilde
K_\lambda(x)=0.
$$
Hence the operator $\widetilde{\Delta_{\lambda}}(\tau)$ can be
reduced to a Hermite operator in $\R^4$:
\begin{equation}\label{hop}
\widetilde{H_\lambda}(\tau)=\sum_{l=1}^{4}\partial_{x_l}^2-4|x|^2\sum_{k=1}^{3}\tau_k^2
-4\sum_{k=1}^{3}\lambda_k\tau_k.
\end{equation}
Now by a result in~\cite{CChT}, we know that the fundamental
solution $\widetilde K_\lambda(x)$ of the operator
$\widetilde{H_\lambda}(\tau)$ has the following form
$$
\widetilde K_\lambda(x,\tau)=\frac{|\tau|^2}{\pi^2}\int_0^\infty
{e^{-4(\sum_{k=1}^3\lambda_k\tau_k)s}}\frac
{e^{-|\tau||x|^2\coth(4|\tau|s)}}{\sinh^2(4|\tau|s)}ds
$$
where $|\tau|^2=\sum_{k=1}^3\tau_k^2$. It follows that
\begin{equation}
\label{eq:klambda}
\begin{split}
K_\lambda(x,t)=&\frac{1}{(2\pi)^3}\int_{\R^3}e^{\sum_{k=1}^3i_kt_k\tau_k}\widetilde
K_\lambda(x,\tau)d\tau\\
&=\frac{1}{8\pi^5}\int_{\R^3}\int_0^\infty\frac{|\tau|^2}{\sinh^2(4|\tau|s)}
e ^{-|\tau||x|^2\coth(4|\tau|s)}e^{\sum_{k=1}^3
(i_kt_k-4\lambda_ks)\tau_k}dsd\tau.
\end{split}
\end{equation}
Let us look at these integrals more carefully. Changing variable
$4|\tau|s=u$ implies that
$$
\widetilde K_\lambda(x,\tau)=\frac{|\tau|}{4\pi^2}\int_0^\infty
{e^{-(\sum_{k=1}^3\lambda_k\tau_k)\frac{u}{|\tau|}}} \frac
{e^{-|\tau||x|^2\coth(u)}}{\sinh^2(u)} du.
$$
Then
\begin{equation*}
\begin{split}
K_\lambda(x,t)=&\frac{1}{(2\pi)^3}\int_{\R^3}e^{\sum_{k=1}^3i_kt_k\tau_k}\widetilde
K_\lambda(x,\tau)d\tau\\
=&\frac{1}{(2\pi)^5}\int_{\R^3}e^{\sum_{k=1}^3i_kt_k\tau_k}
\int_0^\infty {e^{-(\sum_{k=1}^3\lambda_k\tau_k)\frac{u}{|\tau|}}}
\frac {e^{-|\tau||x|^2\coth(u)}}{\sinh^2(u)}|\tau|\,du\,d\tau
\end{split}
\end{equation*}
We introduce the polar coordinates for the $\tau$-variable such that
$\tau=|\tau|\vec n=|\tau|(n_1,n_2,n_3)=r\vec n$ where $r=|\tau|$.
Then the above integral can be rewritten as
\begin{equation*}
\begin{split}
K_\lambda(x,t)=&\frac{1}{(2\pi)^5}\int_0^\infty
\int_{S^2}e^{r(\sum_{k=1}^3i_kt_kn_k)}\int_0^\infty
e^{-(\sum_{k=1}^3\lambda_kn_k)u-r|x|^2\coth(u)}\frac{r^3}{\sinh^2(u)}
\,du\,d\sigma\, dr\\
=&\frac{1}{(2\pi)^5}\int_{S^2}\int_0^\infty
\frac{e^{-(\sum_{k=1}^3\lambda_kn_k)u}}{\sinh^2(u)}\int_0^\infty
e^{-(|x|^2\coth(u)-\sum_{k=1}^3i_kt_kn_k)r}r^3\,dr\, du\,d\sigma\\
=&\frac{6}{(2\pi)^5}\int_{S^2}\int_0^\infty
\frac{e^{-(\sum_{k=1}^3\lambda_kn_k)u}}{\sinh^2(u)}
\frac{du\,d\sigma}{[|x|^2\coth(u)-\sum_{k=1}^3i_kt_kn_k]^4}.
\end{split}
\end{equation*}
Here $d\sigma$ is the surface measure on $S^2$. When $\lambda_k=0$
for $k=1,2,3$, the integral can be simplified. Let $v=\coth(u)$,
then $dv=-\frac {du}{\sinh^2(u)}$ and the above inner integral
reduces to
$$
\int_0^\infty\frac{dv}{[|x|^2v-\sum_{k=1}^3i_kt_kn_k]^4}
=-\frac{1}{3|x|^2}\frac{1}{[|x|^2-\sum_{k=1}^3i_kt_kn_k]^3}.
$$
Hence
\begin{equation}\label{k0}
K_0(x,t)=-\frac{2}{(2\pi)^5|x|^2}\int_{S^2}\frac{d\sigma}
{[|x|^2-\sum_{k=1}^3i_kt_kn_k]^3}.
\end{equation}

If we put $i_1=i$, the usual complex unity and suppose that $i_2,
i_3$ are absent, then the group $\mathcal Q$ is reduced to the
$2$-dimensional Heisenberg group ${\mathbb C}^2\times {\mathbb R}$
with a $4$-dimensional horizontal space and an one dimensional
center. In this case the formula~\eqref{eq:klambda} reduced to the
known formula obtained in~\cite{FS}. Indeed, the Hermit
operator~\eqref{hop} is written in the form
$$\widetilde{H_\lambda}(\tau)=\sum_{l=1}^{4}\partial_{x_l}^2-4|x|^2\tau^2
-4\lambda\tau.
$$ Then
\begin{equation*}
\begin{split}
K_\lambda(x)=&\frac{1}{2\pi}\int_{\R}e^{it\tau}\widetilde
K_\lambda(x)d\tau\\
=&\frac{1}{2\pi^3}\int_{\R}e^{it\tau} \int_0^\infty e^{-4\lambda
\tau s} \frac{|\tau|^2}{\sinh^2(4|\tau| s)}
e^{-|\tau||x|^2\coth(4|\tau|s)}\,ds\,d\tau
\end{split}
\end{equation*}
Changing variable $4s$ to $s$, the above formula becomes
\begin{equation*}
\begin{split}
K_\lambda(x,t)=&\frac{1}{8\pi^3}\int_{\R}e^{it\tau}
|\tau|^2\,\int_0^\infty \frac{e^{-\lambda \tau s}}{\sinh^2(|\tau|
s)}
e^{-|\tau||x|^2\coth(|\tau|s)}\,ds\,d\tau\\
=&\frac{1}{8\pi^3}\int_{\R}e^{it\tau} \tau^2\,\int_0^\infty
\frac{e^{-\lambda \tau s}}{\sinh^2(\tau s)} e^{-\tau|x|^2\coth(\tau
s)}\,ds\,d\tau,
\end{split}
\end{equation*}
since $|\tau|^2$, $\sinh^2(|\tau|s)$, $|\tau|\coth(|\tau|s)$ are all
even functions. Let $u=\tau s$, then one has
\begin{equation*}
\begin{split}
K_\lambda(x,t)=&\frac{1}{8\pi^3}\int_{\R}e^{it\frac{u}{s}}\frac{u^2}{s^2}
\int_0^\infty \frac{e^{-\lambda u}}{\sinh^2(u)} e^{-\frac
{u}{s}|x|^2\coth(u)}\,\frac{ds}{s}\,du\\
=&\frac{1}{8\pi^3}\int_{\R}\frac {u^2}{\sinh^2(u)}e^{-\lambda
u}\int_0^\infty s^{-3}e^{-\frac{1}{s}(u|x|^2\coth(u)-it u)}ds\,du.
\end{split}
\end{equation*}
Once again, changing variable $\frac {1}{s}$ to $s$ in the second
integral, then we have
$$
\int_0^\infty s^{-3}e^{-\frac{1}{s}(u|x|^2\coth(u)-it
u)}ds=\int_0^\infty se^{-s(u|x|^2\coth(u)-it u)}ds.
$$
Now we may apply the identity
$$
\frac{1}{\Gamma (m)}\int_0^\infty s^{m-1}e^{-s A}ds=\frac
{1}{A^m}\qquad {\mbox {for}}\quad \re(A)>0
$$
to the above integral, then we have
$$
K_\lambda(x,t)=\frac{\Gamma(2)}{8\pi^3}\int_{-\infty}^{+\infty}\frac{e^{-\lambda
u}}{\sinh^2(u)}\frac {du}{(|x|^2\coth(u)-it)^2}.
$$
Denote
\[
r=(|x|^4+t^2)^{\frac 1{4}}\quad{\mbox {and}}\quad
e^{-i\phi}=r^{-2}(|x|^2-it)
\]
with $\phi\in \left(-\frac {\pi}{2}, \frac {\pi}{2}\right)$. Using
the identity
\[
\cosh(u+i\phi)=\cosh(u)\cos\phi+i\sinh(u)\sin\phi,
\]
one has
\begin{equation}
\label{isoeq} K_\lambda(x,t)=\frac{1}{8\pi^3}
\int_{-\infty}^\infty\frac {e^{-\lambda
u}}{[r^2\cosh(u+i\phi)]^2}du.
\end{equation}
Changing the contour, the formula (\ref{isoeq}) becomes
\[
K_\lambda(x,t)=\frac{1}{8\pi^3r^{4}} e^{i\lambda \phi}
\int_{-\infty}^\infty\frac {e^{-\lambda u}}{[\cosh(u)]^2}du.
\]
The above integral can be evaluated as follows:
\begin{equation*}
\begin{split}
K_\lambda(x,t)&=\frac 1{4\pi^3}{r^{-4}} e^{i\lambda\phi}
\Gamma\left (\frac {2+\lambda}{2}\right)
\Gamma\left (\frac {2-\lambda}{2}\right)\\
&=\frac {1}{4\pi^3}\Gamma\left (\frac {2+\lambda}{2}\right)
\Gamma\left (\frac {2-\lambda}{2}\right)(|x|^2-it)^{-\frac
{2+\lambda}{2}} (|x|^2+it)^{-\frac {2-\lambda}{2}}.
\end{split}
\end{equation*}
From the above formula, we know that the kernel can be extended from
$|\re(\lambda)|<2$ to the region $\C\setminus \Lambda$ where
\[
\Lambda=\left\{\pm (2+2k):\,\, k\in {\mathbb Z}_+\right\}.
\]
This coincides with the result obtained by Folland and
Stein~\cite{FS}. In particular, when $\lambda=0$, one has
\[
K_0(x,t)=\frac 1{4\pi^3}(|x|^2+it)^{-\frac {2}{2}}
(|x|^2-it)^{-\frac {2}{2}}=\frac 1{4\pi^3}(|x|^4+t^2)^{-1}.
\]
which recovered a result of Folland~\cite{F1}.

\begin{remark}
Let us look on Cauchy-Riemann equations on a domain
$D\subset\mathbb H^2$. Let $u(q_1,q_2)$ be a smooth function of
two quaternion variables, defined in a domain $D\subset{\mathbb
H}^2$. Then the Cauchy-Riemann equations is expressed by the
system
\begin{equation}\label{CRsys}\begin{split}
& \bar\partial_{q_1}u  =\frac{\partial u}{\partial \bar q_1} =f_1
\\& \bar\partial_{q_2}u =\frac{\partial u}{\partial \bar q_2}=f_2.
\end{split}\end{equation}

Since the product $ \bar\partial_{q_1} \bar\partial_{q_2}$ is not
commutative, we can not obtain the compatibility condition on $f$
of the type
$$\frac{\partial f_2}{\partial \bar q_1}=\frac{\partial f_1}{\partial \bar q_2}.$$
So the result $\bar \partial_b^2\neq 0$ is quite different from
the Heisenberg group. Moreover, calculations of forms for
quaternion is more complicated. The simplest anticommutative
property
\begin{equation}\label{anticpr}
dq_1\wedge dq_2=-dq_2\wedge dq_1
\end{equation}
does not hold. This implies that $dq\wedge dq\neq 0$. So even if
$[\bar H,\bar H]=0$ we can not expect $\bar\partial_b^2$,
since~\eqref{anticpr}. It was proved in~\cite{ABLSS} that the
system~\eqref{CRsys} has a $C^{\infty}$ solution $u(q_1,q_2)$
defined in a convex domain $D\subset {\mathbb H}^2$ if and only if
the vector $(f_1,f_2)$ solve a matrix the second order
differential equation that can be considered as an analogue of the
compatibility condition. However, can we find an elegant way to
describe is the operator
$\square_b=\bar\partial_b\bar\partial_b^\ast+\bar\partial_b^\ast\bar\partial_b$?
Furthermore, what is the correct setting for the
$\bar\partial$-Neumann problem on a bounded domain in ${\mathbb
H}^2$?
\end{remark}
\smallskip

\begin{remark} How to define correctly the pseudoconvexity in quaternion setting?
We need the analogue of Hermitian matrix and Levi matrix. Even for
the case of 1-dimension it is need to understand. We have the
tangential Cauchy-Riemann-Fueter operator $\bar H$, such that $[H,\bar
H]=2\sum_{k=1}^{3}i_k\partial_{t_k}$. We can say that
$\sum_{k=1}^{3}i_k\partial_{t_k}$ is an elementary pure imaginary
vector field defining the missing directions of real dimension
$3$. Then the Levi  matrix is the number $2$, which is positive
definite. What is in general case?
\end{remark}

\end{document}